\documentclass[10pt]{article}

\usepackage{amsmath}
\usepackage{amssymb}
\usepackage{graphicx}
\usepackage{epic,eepic}

\usepackage{color}

\usepackage{amstext,amsthm,amssymb,amsmath}
\usepackage{graphicx}
\usepackage{subfigure}
\usepackage{wrapfig}
\usepackage{fullpage}
\usepackage{color}
\usepackage{multirow}
\usepackage{tabulary}
\usepackage{booktabs}
\usepackage{enumerate}

\newcommand{\average}[1]{\{\!\!\!\{{#1}\}\!\!\!\}}
\newcommand{\jump}[1]{[\![{#1}]\!]}

\newtheorem{theorem}{Theorem}[section]
\newtheorem{lemma}[theorem]{Lemma}

\hfuzz=\maxdimen
\title{An online generalized multiscale discontinuous Galerkin method (GMsDGM) for flows in heterogeneous media}

\author{Eric T. Chung\thanks{Department of Mathematics, The Chinese University of Hong Kong, Hong Kong SAR.
This research is partially supported by the Hong Kong RGC General Research Fund (Project number: 400813).}, \;
Yalchin Efendiev\thanks{Department of Mathematics, Texas A\&M University, College Station, TX; Numerical Porous Media SRI Center, King Abdullah University of Science and Technology (KAUST), Thuwal 23955-6900, Kingdom of Saudi Arabia} \;
and \; Wing Tat Leung\thanks{Department of Mathematics, Texas A\&M University, College Station, TX.}
}

\begin{document}
\maketitle

\begin{abstract}

Offline computation is an essential component in
 most multiscale model reduction techniques.
However, there are multiscale
problems in which offline procedure is insufficient to
give accurate representations of solutions,
due to the fact that offline computations are typically performed
locally and global information is missing in these offline information.
To tackle this difficulty,
we develop an online local adaptivity technique
for local multiscale model reduction problems.
We design new online basis functions within Discontinuous Galerkin
method based on local residuals and
some optimally estimates.
The resulting basis functions are able to capture the solution efficiently
and accurately, and are added to the approximation iteratively.
Moreover, we show that the iterative procedure is convergent with a rate independent of physical scales if the initial space is chosen carefully.
Our analysis also gives a guideline on how to choose the initial space.
We  present some numerical examples to show the performance
of the proposed method.

\end{abstract}

\section{Introduction}
\label{sec:intro}

In this paper, we develop an online local adaptivity technique
for a class of multiscale model reduction problems.
Many realistic applications involve solving problems that contain multiple scales and high contrast.
Direct solution methods for these problems require fine-grid discretizations and result in large discrete systems
that are computationally intractable.
Common model reduction techniques perform
the discretization of the problems on a coarse
grid, which is much larger than the scales under consideration,
with the aim of getting more efficient solution strategies.
There are a variety of multiscale model reduction techniques based on numerical upscaling
(e.g., \cite{dur91, weh02}) or multiscale methods
(e.g., \cite{Arbogast_two_scale_04, Chu_Hou_MathComp_10,ee03,
  egw10,eh09,ehg04, GhommemJCP2013,ReducedCon,MsDG,Wave,WaveGMsFEM}).
Most of the existing techniques are based on the so called offline construction.
In particular, reduced models are computed in a pre-processing step, called offline stage,
before the actual simulations, called online stage, are performed.
For instances, some effective media are pre-computed for methods based on numerical upscaling
and some multiscale basis functions are pre-computed for multiscale finite element methods.
While these methods are effective in a wide variety of applications,
there are still situations for which these methods are inadequate to give reliable solutions unless a large dimensional offline space is employed.
Some of these situations involve external source effects
and distant effects, which are ignored by most multiscale model reduction methods since they are typically based on local constructions.
Therefore, it is evident that offline procedures are sometimes not enough to give efficient
reduced models.
Hence, it is the purpose of this paper to design a novel multiscale model reduction method.
Our proposed method is based on a combination of offline technique and
an online enrichment technique.
The online technique is able to produce a reduced model taking care of external sources and distant effects,
without using global models.
The online construction is also performed locally and adaptively in regions with more heterogeneities, giving very efficient reduced models.

%In this paper, we develop an adaptive Generalized Multiscale
%Discontinuous Galerkin Method (GMsDGM)
%for a class of high-contrast flow problems, and
%derive a-priori and a-posteriori error estimates for the method.
%We propose an adaptive enrichment algorithm for our GMsDGM
%based on the a-posteriori error estimator
%and prove its convergence.
%The enrichment is done using inexpensive $L_2$-based error
%indicators which allows adding
%more basis functions in an automatic way.

Our proposed method follows the overall idea of the
Generalized Multiscale Finite Element Method (GMsFEM), which is introduced in \cite{egh12} and
is a generalization of the classical multiscale finite element method
(\cite{hw97}) in the way that the coarse spaces are systematically enriched, taking into account
small scale information and
complex input spaces.
Instead of conforming finite element spaces as in \cite{egh12,hw97}, we will use in this paper discontinuous Galerkin finite element
spaces, which have some essential advantages (see \cite{eglmsMSDG,WaveGMsFEM}) in multiscale simulations because it allows coupling discontinuous basis
functions.
The discretization starts with a coarse grid and a space of snapshot functions,
which are defined on coarse elements.
A space reduction is then performed to obtain a much smaller offline space
by means of spectral decomposition.
The spectral decomposition is performed locally on coarse elements,
thus the functions in the offline space are in general discontinuous across coarse edges.
The offline space is used as the approximation space for the interior penalty discontinuous Galerkin (IPDG) discretization
on the coarse grid for the problem under consideration, giving our generalized multiscale discontinuous Galerkin method (GMsDGM).
We remark that the offline space is computed only once in the pre-processing offline stage,
and the same set of basis functions is used for any given source terms and boundary conditions.
A-priori error estimate can be derived as in \cite{egw10, eglp13, WaveGMsFEM}
showing that the error
is inverse proportional to the first eigenvalue corresponding
to the first eigenfunction that is not used
in the construction of the reduced space.
Since the aim of the paper is the new online locally adaptive procedure and its convergence,
we will not discuss a-priori error estimate in this paper.

The previous paragraph discusses the offline component of our method.
As we discussed before, some new basis functions are necessary to capture certain behavior
of the solution which cannot be captured by offline basis functions.
For example, the solution may contain heterogeneities due to some distant effects and source terms, and these cannot be incorporated efficiently
by offline basis functions before the solution is computed.
Hence, it is the purpose of this paper to develop a technique to find
new basis functions in the online stage.
Our method consists of an iterative procedure.
Given an approximate solution, some local residuals on coarse elements can be computed
to reflect the amount of error in these coarse elements.
These local residuals serve as indicators to locate regions,
 where new online basis functions are necessary.
We will show that the projection of these residuals to the fine-grid can be used as new basis functions
and that the energy-norm error has the most decay in a certain sense when these residual-based basis functions are included in the next solution process.
In addition, we will show that this iterative procedure is convergent
with a convergence rate independent of scales and contrast.
In our analysis of convergence, we will show that it is essential to choose the appropriate space
to begin the iterative procedure (cf. \cite{ge09_2, ge09_1reduceddim, Efendiev_GLW_ESAIM_12}).
This initial space is computed in the offline stage and is obtained from a carefully design spectral problem.
With this choice of the initial space, we show that one can obtain a very fast decay of errors
by adding our online basis functions.
Finally, we remark that there are offline adaptive enrichment strategies in the context of GMsFEM.
In particular, in \cite{Adaptive-GMsFEM}, offline adaptive procedure is developed and its convergence is analyzed
using techniques in \cite{BrennerScott,AdaptiveFEM}.
This is an efficient method to adaptively enrich the offline space and is desirable for problems,
 where offline basis functions are good enough to capture the solution.
On the other hand, we remark that other adaptive methods are available \cite{Dorfler96,ohl12, abdul_yun, dinh13, nguyen13, tonn11, donoho_sparsity_review}.
Also, we remark that online basis functions within continuous Galerkin GMsFEM
is studied in \cite{cel15_1}.

The rest of the paper is organized in the following way.
In the next section, we present the basic idea of GMsDGM and our online locally adaptive procedure.
The method is then detailed and analyzed
in Section \ref{sec:online}.
In Section \ref{sec:numerresults}, numerical
results are illustrated to test the performance of this adaptive algorithm.
Finally a conclusion is given in Section \ref{sec:conclusion}.

%==============================
\section{Preliminaries}

\label{prelim} %==============================

In this paper,
we consider the following
high-contrast flow problem
\begin{equation}
-\mbox{div}\big(\kappa(x)\,\nabla u\big)=f\quad\text{in}\quad D,\label{eq:original}
\end{equation}
subject to the homogeneous Dirichlet boundary condition $u=g$ on
$\partial D$, where $D$ is the computational domain and $f(x)$ is a given source term.
We assume that the coefficient $\kappa(x)$ is highly heterogeneous with very high contrast.
For the convenience of our analysis, we also assume that $\kappa(x)$ is bounded below, that is, $\kappa(x) \geq 1$.
Due to the heterogeneity and high contrast of $\kappa(x)$, very fine meshes are necessary to obtain accurate numerical solutions.
It is therefore crucial to develop a numerical scheme with a low dimensional approximation space
for the efficient approximation of (\ref{eq:original}).

\begin{figure}[htb]
  \centering
  \includegraphics[width=0.65 \textwidth]{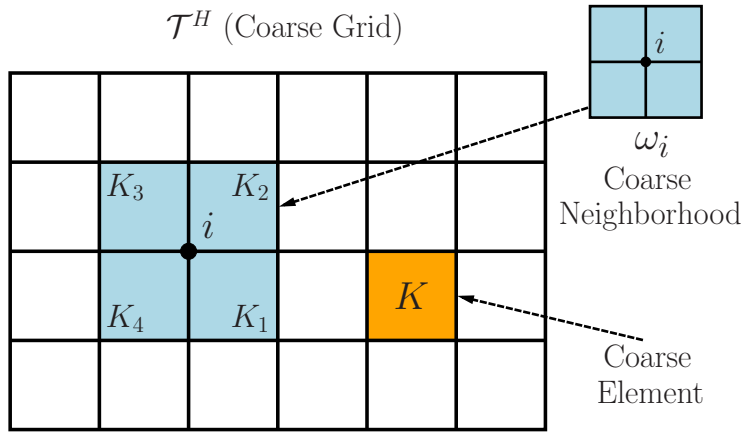}
  \caption{Illustration of a coarse neighborhood and a coarse element.}
  \label{schematic}
\end{figure}

Next, we present some notations needed for the constructions of our scheme.
Consider a given triangulation $\mathcal{T}^{H}$ of the domain $D$
with mesh size $H>0$.
For convenience, we assume the domain $D$ is rectangular and that the triangulation $\mathcal{T}^H$
consists of rectangles.
We call $\mathcal{T}^{H}$ the coarse
grid and $H$ the coarse mesh size. Elements of $\mathcal{T}^{H}$
are called coarse grid blocks and we use $N$ to denote the number of coarse grid blocks.
The set of all coarse grid edges is
denoted by $\mathcal{E}^{H}$.
See Figure~\ref{schematic} for an illustration.
%and the set of all coarse grid nodes
%is denoted by $\mathcal{S}^{H}$. We also use $N_{S}$ to denote the
%number of coarse grid nodes, $N_{E}$ to denote the number of coarse
%grid blocks and $N$ to denote the number of coarse grid blocks.
We also introduce a finer triangulation $\mathcal{T}^{h}$ of the computational domain $D$,
obtained by
a conforming refinement of
the coarse grid $\mathcal{T}^{H}$. We call $\mathcal{T}^{h}$ the
fine grid and $h>0$ the fine mesh size.
%We remark that the use of
%the conforming refinement is only to simplify the discussion of the
%methodology and is not a restriction of the method.

Now we present the framework of our GMsDGM.
%with a residual based locally online adaptivity.
The methodology consists
of two main ingredients, namely, the construction of local basis functions
and the global coarse grid level coupling.
For the coarse grid level coupling, we will apply the interior penalty discontinuous Galerkin (IPDG) method \cite{IPDGbook}.
We remark that other discretizations can also be used.
Assume that $V_H$ is a given approximation space defined on the coarse grid $\mathcal{T}^H$.
Functions in $V_H$ are piecewise polynomials within coarse grid blocks and are in general discontinuous across coarse grid edges.
Following standard procedures, the IPDG method reads:
find $u_{H}\in V_H$
such that
\begin{equation}
a_{\text{DG}}(u_{H},v)=(f,v),\quad\forall v\in V_H,\label{eq:ipdg}
\end{equation}
where the bilinear form $a_{\text{DG}}$ is defined as
\begin{equation}
a_{\text{DG}}(u,v)=a_{H}(u,v)-\sum_{E\in\mathcal{E}^{H}}\int_{E}\Big(\average{{\kappa}\nabla{u}\cdot{n}_{E}}\jump{v}+\average{{\kappa}\nabla{v}\cdot{n}_{E}}\jump{u}\Big)+\sum_{E\in\mathcal{E}^{H}}\frac{\gamma}{h}\int_{E}\overline{\kappa}\jump{u} \jump{v} \label{eq:bilinear-ipdg}
\end{equation}
with
\begin{equation}
a_{H}({u},{v})=\sum_{K\in\mathcal{T}_{H}}a_{H}^{K}(u,v),\quad a_{H}^{K}(u,v)=\int_{K}\kappa\nabla u\cdot\nabla v,
\end{equation}
where $\gamma>0$ is a penalty parameter, ${n}_{E}$ is a fixed unit
normal vector defined on the coarse edge $E \in \mathcal{E}^H$.
Note that, in (\ref{eq:bilinear-ipdg}),
the average and the jump operators are defined in the classical way.
Specifically, consider an interior coarse edge $E\in\mathcal{E}^{H}$
and let $K^{+}$ and $K^{-}$ be the two coarse grid blocks sharing
the edge $E$. For a piecewise smooth function $G$ with respect to the coarse grid $\mathcal{T}^H$, we define
\[
\average{G}=\frac{1}{2}(G^{+}+G^{-}),\quad\quad\jump{G}=G^{+}-G^{-},\quad\quad\text{ on }\, E,
\]
where $G^{+}=G|_{K^{+}}$ and $G^{-}=G|_{K^{-}}$ and we assume that
the normal vector ${n}_{E}$ is pointing from $K^{+}$ to $K^{-}$.
Moreover, on the edge $E$, we define $\overline{\kappa} = (\kappa_{K^+}+\kappa_{K^-})/2$,
where $\kappa_{K^{\pm}}$ is the maximum value of $\kappa$ over $K^{\pm}$.
For a coarse edge $E$ lying on the boundary $\partial D$, we define
\[
\average{G}=\jump{G}=G,\quad \text{ and }\quad \overline{\kappa} = \kappa_{K} \quad\quad\text{ on }\, E,
\]
where we always assume that ${n}_{E}$ is pointing outside of $D$.

For our analysis, we define the DG-norm as
\begin{equation*}
\| u\|_{\text{DG}}^2 =
a_{H}({u},{u})+\sum_{E\in\mathcal{E}^{H}}\frac{\gamma}{h}\int_{E}\overline{\kappa}\jump{u}^2.
\end{equation*}
Then, the following continuity and coercivity of the bilinear form $a_{\text{DG}}$ hold.
For completeness, we include a proof of this result in the Appendix.
%The proof can be found in the Appendix.

\begin{lemma}
\label{lem:norm}
Assume that the penalty parameter $\gamma$ is chosen
so that $\gamma>C_{\text{\rm inv}}^{2}$. The bilinear form $a_{\text{DG}}$
defined in (\ref{eq:bilinear-ipdg}) is continuous and coercive with respect to the DG-norm, that
is,
\begin{eqnarray}
a_{\text{DG}}({u},{v}) & \leq & a_{1}\|{u}\|_{\text{DG}}\,\|{v}\|_{\text{DG}}, \label{lem:cont} \\
a_{\text{DG}}({u},{u}) & \geq & a_{0}\|{u}\|_{\text{DG}}^{2}, \label{lem:coer}
\end{eqnarray}
for all ${u},{v} \in V_H$, where $a_{0}=1-C_{\text{\rm inv}}\gamma^{-\frac{1}{2}}>0$ and $a_{1}=1+C_{\text{\rm inv}}\gamma^{-\frac{1}{2}}$.
\end{lemma}

One main result of the paper is a convergence estimate of an adaptive procedure for the problem (\ref{eq:ipdg}).
For this purpose, we will compare the multiscale solution $u_H$ to a fine-scale solution $u_h$ defined in the following way.
We first let
\begin{equation*}
V^h_{\text{DG}} = \{ v \in L^2(D) \, : \, v|_{K} \in V^h(K) \},
\end{equation*}
where $V^h(K)$ is the space of continuous piecewise bilinear functions defined on $K$ with respect to the fine grid.
The fine-scale solution $u_{h}\in V_{\text{DG}}^{h}$ is defined as the solution of the following
\begin{equation}
a_{\text{DG}}(u_{h},v)=(f,v),\quad\forall v\in V_{\text{DG}}^{h}.\label{eq:ipdgfine}
\end{equation}
It is well-known that $u_{h}$ gives a good approximation to the exact solution $u$
up to a coarse grid discretization error.
%in the DG-norm as the fine mesh size $h\rightarrow0$.

The second main component of our method is
the construction of local basis functions, which contains two stages, namely the offline stage and the online stage.
In the offline stage,
a snapshot space $V^{i,\text{snap}}$ is
first constructed for each coarse grid block $K_{i} \in\mathcal{T}^H$.
The snapshot space contains a rich space of basis functions, which can be used to approximate
the fine-scale solution defined (\ref{eq:ipdgfine}) with a good accuracy.
A spectral problem
is then solved in the snapshot space $V^{i,\text{snap}}$ and eigenfunctions corresponding
to dominant modes are used as the basis functions. The resulting space
is called the local offline space $V^{i,\text{off}}$ for the $i$-th coarse
grid block $K_i$. The global offline space $V^{\text{off}}$ is then defined
as the linear span of all these $V^{i,\text{off}}$, for $i=1,2,\cdots,N$.
This global offline space $V^{\text{off}}$ will be used as the initial space of
our method.
%which can be
%formulated as:
%For vector-valued functions, the above average and jump operators
%are defined component-wise.
%We note that the DG coupling (\ref{eq:ipdg})
%is the classical interior penalty discontinuous Galerkin (IPDG) method \cite{IPDGbook}
%with our multiscale basis functions as the approximation space.
We denote this initial space as $V_H^{(0)}$.
Using the initial space, an initial solution $u_H^{(0)}$ can be computed by solving (\ref{eq:ipdg}).
Local residuals in coarse grid blocks can then be computed based on the initial solution $u_H^{(0)}$.
In coarse grid blocks with large residuals, new basis functions are computed
and added to the approximation space. This procedure is continued until certain tolerance is reached.
Next, we present a general outline of the method.

Assume that the initial space $V_H^{(0)}$ is given
and the initial solution $u_H^{(0)}$ is computed.
For any $m \geq 0$, we repeat the following until the solution $u_H^{(m)}$ satisfies certain tolerance requirement.
\begin{enumerate}
\item[Step 1:] Solve (\ref{eq:ipdg}) using the space $V_H^{(m)}$ to obtain the solution $u_H^{(m)} \in V_H^{(m)}$.
\item[Step 2:] Compute local residuals based on the solution $u_H^{(m)}$.
\item[Step 3:] Construct new basis functions in regions, where the residuals are large.
\item[Step 4:] Add these basis functions to $V_H^{(m)}$ to form a new space $V^{(m+1)}$.
\end{enumerate}

In the following, we will give the details of Step 2 and Step 3.
We will also explain how one chooses the initial space $V_H^{(0)}$.

%============================================

\section{Locally online adaptivity}
\label{sec:online}
In this section, we will give details of our locally online adaptivity for the problem (\ref{eq:ipdg}).
As presented in the general outline of the method from the previous section,
our adaptivity idea contains the choice of initial space
as well as construction of new local multiscale basis functions.
In the following, we will give the construction of these in detail.

\subsection{Initial space}

We present the definition of the initial space $V_H^{(0)}$.
Let $x_i$ be a node in the coarse grid $\mathcal{T}^H$, referred to as the $i$-th coarse node,
for $i=1,2,\cdots, N_c$,
where $N_c$ is the number of nodes in the coarse grid $\mathcal{T}^H$.
We will then define the $i$-th coarse neighbourhood $\omega_i$ as the union of all coarse grid blocks
having the node $x_i$, see Figure~\ref{schematic}.
Moreover, for each coarse grid block $K\in\mathcal{T}^H$, we let $\chi^K_{(j)}$, $j=1,2,3,4$,
be the partition of unity functions, having value $1$ at one vertex $y_j$ and value $0$ at the remaining three vertices,
where $y_j$, $j=1,2,3,4$, are the four vertices of $K$.
Note that
there is exactly one value of $j$ such that the vertex $y_j$ is the same as the vertex $x_i$.
In the case, we write $\chi_{(j)}^K = \chi_i^K$.
One can use the standard multiscale basis functions or bilinear functions
as the partition of unity functions.
Note that we do not require any continuity of these partition of unity functions across coarse grid edges.
The partition of unity functions are all supported on coarse grid blocks.
Furthermore, we define the space $V^h(\omega_i)$ by
\begin{equation*}
V^h(\omega_i) = \{ v \in L^2(\omega_i) \, : \, v|_{K} \in V^h(K), \; K\in\mathcal{T}^H, \; K\subset\omega_i \}.
\end{equation*}
That is, functions in $V^h(\omega_i)$ are supported in $\omega_i$ and belong to the space $V^h(K)$
for each coarse grid block $K\subset \omega_i$.
Note that there is no continuity condition across boundaries of coarse grid blocks.
We consider $V^h(\omega_i)$ as the snapshot space in $\omega_i$, that is $V^{i,\text{snap}} = V^h(\omega_i)$,
and perform a dimension reduction through a spectral problem.
For this purpose, we define $\mathcal{E}^H_i$ be the set of coarse grid edges lying in the interior of $\omega_i$, and
the following bilinear form
\begin{equation}
a_{\omega_i}(u,v) = \sum_{K\in\mathcal{T}^H,K\subset \omega_i} a_H^K(u,v) + \sum_{E\in \mathcal{E}^H_i} \frac{\gamma}{h} \int_E \overline{\kappa} \jump{u} \jump{v}, \quad
\forall u, v\in V^h(\omega_i).
\end{equation}
Based on our analysis to be presented next, we solve the following spectral problem
\begin{equation}
a_{\omega_i}(u,v) = \lambda s_{\omega_i}(u,v), \quad \forall v\in V^h(\omega_i),
\label{eq:spec}
\end{equation}
where
\begin{equation}
s_{\omega_i}(u,v) = \sum_{K\in\mathcal{T}^H,K\subset \omega_i} \int_{K} \kappa  |\nabla \chi^K_i|^2  u \, v
+ \sum_{E\in\mathcal{E}^H_i} \frac{\gamma}{h} \int_E \overline{\kappa}   \jump{\chi_i^K}^2  \average{u} \, \average{v}, \quad \forall u,v\in V^h(\omega_i).
\end{equation}
We use the notations $\lambda_k^{\omega_i}$ and $\Psi_k^{\omega_i}$
to denote the $k$-th eigenvalue and the $k$-th eigenvector of the above spectral problem (\ref{eq:spec}).
Each eigenfunction $\Psi_k^{\omega_i}$ corresponds to a function in $\psi_k^{\omega_i} \in V^h(\omega_i)$ defined by
\begin{equation*}
\psi_k^{\omega_i} = \sum_{j=1}^{n_i} (\Psi_k^{\omega_i})_j w_j^{\omega_i},
\end{equation*}
where $n_i$ is the dimension of $V^h(\omega_i)$ and $\{ w_j^{\omega_i} \}_{j=1}^{n_i}$
is a basis for $V^h(\omega_i)$.
In the above definition, $(\Psi_k^{\omega_i})_j$ is the $j$-th component of the eigenvector $\Psi_k^{\omega_i}$.

\begin{figure}[htb]
  \centering
  \includegraphics[width=3in, height=3in]{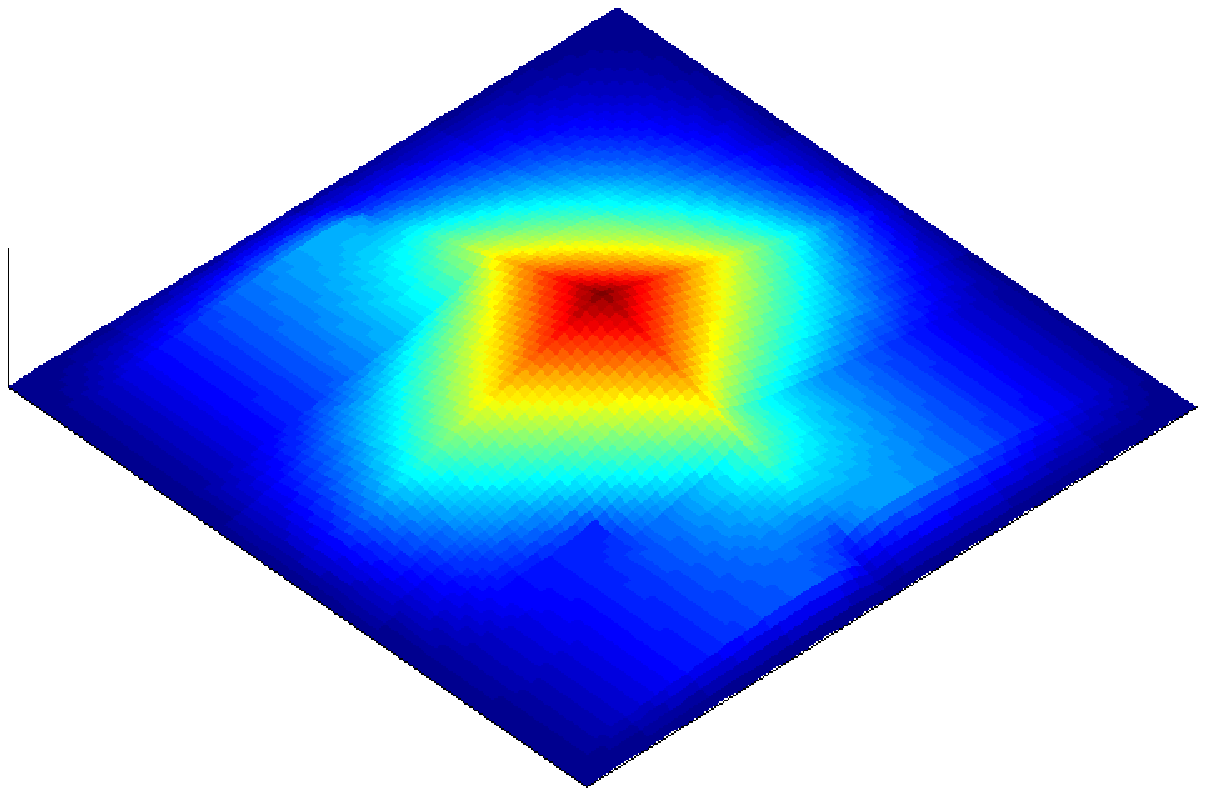}
  \includegraphics[width=3in, height=3in]{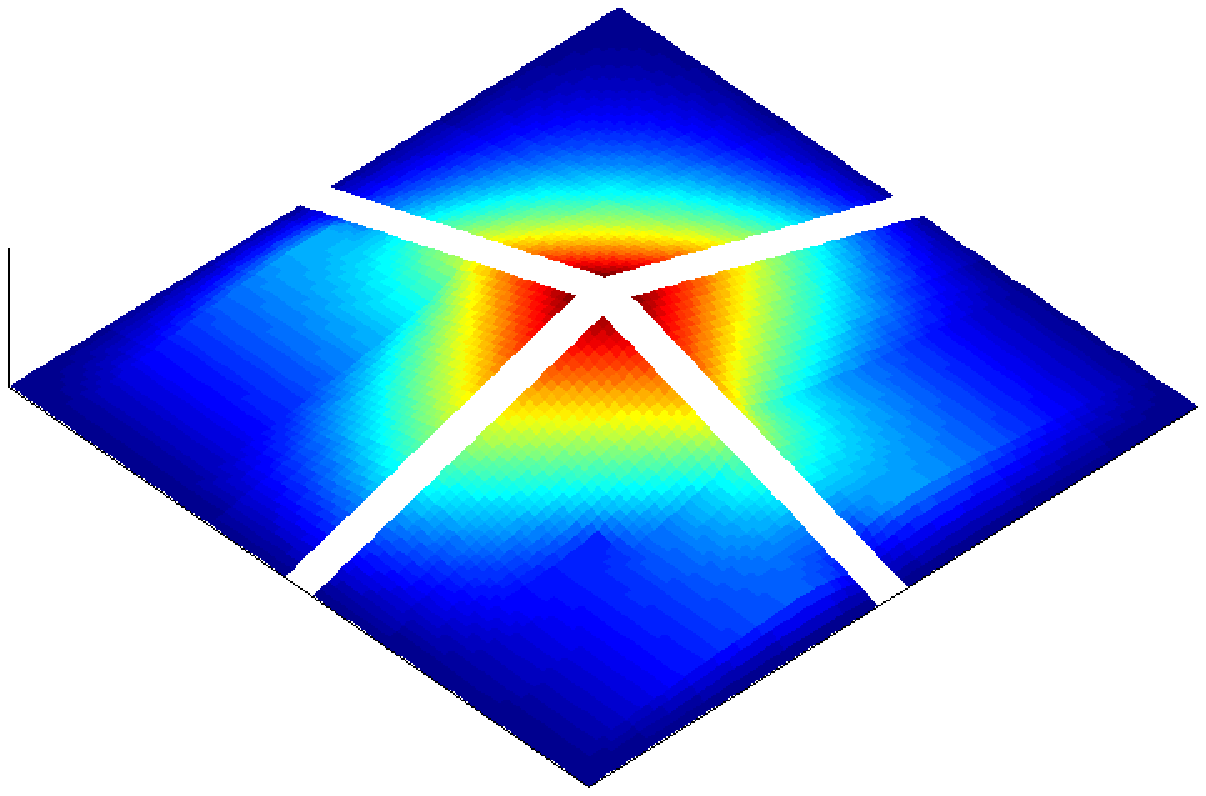}
  \caption{Illustration of the initial basis construction. Left: An eigenfunction $\chi_i^K \psi_k^{\omega_i}$ is defined in $\omega_i$. Right: $4$ basis functions are obtained by splitting $\chi_i^K \psi_k^{\omega_i}$
   into $4$ pieces, and each has support in $K\subset \omega_i$.}
  \label{schematic1}
\end{figure}

For each $\omega_i$, we solve the spectral problem (\ref{eq:spec})
and the first $L_i$ eigenfunctions are used to form the initial space.
Each eigenfunction $\psi_k^{\omega_i}$ will be first multiplied by the partition of unity function $\chi_i^K$, for each $K\subset\omega_i$,
and is then decoupled across coarse grid edges to form $4$ basis functions.
In particular, the $4$ new basis functions have support in one of the coarse grid block forming $\omega_i$
and are zero in the other three coarse grid blocks forming $\omega_i$.
For example, if $K\subset \omega_i$, the basis function is $\chi_i^K \psi_k^{\omega_i}$.
We write $V^{i,\text{off}}$ as the space spanned by $\chi_i^K \psi_k^{\omega_i}$, for all $K\subset\omega_i$.
See Figure~\ref{schematic1} for an illustration.
The initial space $V_H^{(0)}$ is obtained by the linear span of all functions constructed in the above procedure.
%(see Figure \ref{schematic1} for illustration).

\subsection{Construction of online basis functions}

In this section, we will discuss the construction of our local online basis functions.
The purpose is to add basis functions locally in some coarse neighborhoods
to obtain rapidly decaying errors.
Assume that the space $V_H^{(m)}$ at the $m$-th iteration and the corresponding solution $u_H^{(m)}$ are given.
For each coarse neighborhood $\omega_i$, we define the local residual by
\begin{equation}
R^{(m)}_i(v) =  (f, v) - a_{\text{DG}}(u_H^{(m)}, v), \quad v \in V_0^h(\omega_i)
\label{eq:res}
\end{equation}
where $V_0^h(\omega_i) \subset V^h(\omega_i)$ contains functions that are zero on $\partial\omega_i$.
The residual $R_i^{(m)}$ can be seen as a linear functional defined on $V^h_0(\omega_i)$
with norm $\| R_i^{(m)} \|$ defined by
\begin{equation*}
\| R_i^{(m)} \| = \sup_{v\in V_0^h(\omega_i)} \frac{| R_i^{(m)}(v)|}{\|  v\|_{\omega_i}},
\end{equation*}
and $\|v\|_{\omega_i}^2 = a_{\omega_i}(v,v)$.
%We also define the restriction of the bilinear form $a_{\text{DG}}$ to $V^h(\omega_i)$ as $a_{\text{DG},\omega_i}$.
%Note that $a_{\text{DG},\omega_i}$ and $a_{\omega_i}$ are not the same,
%since $a_{\text{DG},\omega_i}$ contains integrals on edges lying on the boundary of $\omega_i$.
We will then find the new online basis function $\phi \in V_0^h(\omega_i)$ by solving
\begin{equation}
 a_{\omega_i}( \phi, v) = R_i^{(m)}(v), \quad \forall v\in V_0^h(\omega_i).
\label{eq:online}
\end{equation}
%In the above definition (\ref{eq:online}), $\chi_i^K$ is considered to be defined only on $K$, and has zero value outside $K$.
%Next, we define the following elliptic projection operator $P_{\text{DG}}: H^1(\mathcal{T}^H) \rightarrow V^h_{\text{DG}}$ by
%\begin{equation*}
%a_{\text{DG}}(P_{\text{DG}}w,v) = a_{\text{DG}}(w,v), \quad\forall v\in V^h_{\text{DG}},
%\end{equation*}
%where $H^1(\mathcal{T}^H)$ is the space of piecewise $H^1$ functions with respect to the coarse grid $\mathcal{T}^H$.
%It is then easy to show that $\|P_{\text{DG}} w\|_A \leq \|w\|_A$ for all $w\in H^1(\mathcal{T}^H)$,
%where the A-norm is defined by
%\begin{equation*}
%\| u\|_A^2 = a_{\text{DG}}(u,u), \quad \forall u\in V_H.
%\end{equation*}
The new online basis function $\phi$ is added to $V_H^{(m)}$ to form $V_H^{(m+1)}$.
%Note that $\widetilde{\phi} \in V^h(\omega_i)$.

%\subsection{Convergence}

The motivation of finding a new basis function $\phi$ by solving (\ref{eq:online})
can be explained as follows.
We define the A-norm by
\begin{equation*}
\| u\|_A^2 = a_{\text{DG}}(u,u), \quad \forall u\in V_H.
\end{equation*}
We notice that
the A-norm is equivalent to the DG-norm $\|u\|_{\text{DG}}$ by Lemma \ref{lem:norm}.
From (\ref{eq:ipdg}) and (\ref{eq:ipdgfine}), we have
the following Galerkin orthogonality condition
\begin{equation}
a_{\text{DG}}(u_h - u^{(m)}_H, v) = 0, \quad\forall v\in V^{(m)}_H.
\label{eq:ortho}
\end{equation}
Thus, we see that the following optimal error bound holds
\begin{equation}
\| u_h - u^{(m)}_H \|^2_A \leq  \| u_h - \widetilde{u} \|_A^2, \quad \forall \widetilde{u} \in V^{(m)}_H.
\label{eq:online1}
\end{equation}
%Notice that this error bound holds if we replace $V_H$ by any
%$V_H^{(m)}$, $m\geq 0$.
Notice that (\ref{eq:ortho}) and (\ref{eq:online1}) hold for any $m \geq 0$.

%We let $V^{(m)}$ be the approximation space in the $m$-th step.
%We will next construction the approximation space $V^{(m+1)}$ for the $(m+1)$-th step.
%Let $\omega_i$ be the coarse neighborhood corresponding to the $i$-th coarse node.
We will enrich the space $V_H^{(m)}$ by adding a basis function $\phi$ in the space $V_0^h(\omega_i)$
to form $V_H^{(m+1)}$.
%and is zero on $\partial\omega_i$.
First, (\ref{eq:online1}) implies
\begin{equation}
\| u_h - u^{(m+1)}_H \|^2_A \leq  \| u_h - \widetilde{u} \|_A^2, \quad \forall \widetilde{u} \in V_H^{(m+1)}.
\label{eq:online2}
\end{equation}
%Assume that $\widetilde{\phi}$ is added in
Taking $\widetilde{u} = u_H^{(m)} + \alpha \phi$, for some scalar $\alpha$, in (\ref{eq:online2}),
we have
\begin{equation*}
\| u_h - u_H^{(m+1)} \|_A^2 \leq  \| u_h - u_H^{(m)} - \alpha \phi\|_A^2
\end{equation*}
which implies
\begin{equation*}
\| u_h - u_H^{(m+1)} \|_A^2 \leq  \| u_h - u_H^{(m)}\|_A^2 - 2\alpha \, a_{\text{DG}}(u_h-u_H^{(m)},\phi) + \alpha^2 \| \phi\|_A^2.
\end{equation*}
%Therefore, by the definition of $P_{\text{DG}}$ and the
%fact that $\|P_{\text{DG}} w\|_A \leq \|w\|_A$, we have
%\begin{equation*}
%\| u_h - u_H^{(m+1)} \|_A^2 \leq  \| u_h - u_H^{(m)}\|_A^2 - 2 \alpha \sum_{K\subset\omega_i}  a_{\text{DG}}(u_h-u_H^{(m)}, \chi_i^K \phi) + \alpha^2 \| \sum_{K\subset\omega_i} \chi_i^K \phi \|_A^2.
%\end{equation*}
%We normalize $\phi$ so that $\|  \phi \|_A=1$.
Taking $\alpha =  a_{\text{DG}}(u_h-u_H^{(m)},  \phi) / \|\phi\|^2_A$, we obtain
\begin{equation}
\| u_h - u_H^{(m+1)} \|_A^2 \leq  \| u_h - u_H^{(m)}\|_A^2 -  \frac{ a_{\text{DG}}(u_h-u_H^{(m)}, \phi)^2}{\|\phi\|_A^2}.
\label{eq:su_error}
\end{equation}
By the definition of the residual $R_i^{(m)}$ in (\ref{eq:res}),
%Given $\omega_i$, we define the local residual by
%\begin{equation*}
%R^{(m)}_i(\phi) = (f,\phi) - a_{\text{DG}}(u_H^{(m)},\phi).
%\end{equation*}
we see that (\ref{eq:su_error}) becomes
\begin{equation}
\| u_h - u_H^{(m+1)} \|_A^2 \leq  \| u_h - u_H^{(m)}\|_A^2 - \frac{(R^{(m)}_i(\phi))^2}{\|\phi\|_A^2}.
\label{eq:su_error1}
\end{equation}
From (\ref{eq:su_error1}), we see that the quantity $(R^{(m)}_i(\phi))^2 / \|\phi\|_A^2$ measures the amount of reduction in error
when the basis function $\phi$ is added in $V_H^{(m)}$ to form $V_H^{(m+1)}$.
We will construct the function $\phi \in V_0^h(\omega_i)$ to obtain the most reduction in error. Thus, we find $\phi \in V_0^h(\omega_i)$
that maximizes $R^{(m)}_i(\phi) / \|\phi\|_A$. Equivalently, we find $\phi \in V_0^h(\omega_i)$ by solving
\begin{equation*}
a_{\omega_i}( \phi, v) = R^{(m)}_i(v), \quad \forall v \in V_0^h(\omega_i).
\end{equation*}
%This problem is the same as (\ref{eq:online}).
Notice that, we have used the fact that $\|\phi\|_A = \|\phi\|_{\omega_i}$ when $\phi\in V^h_0(\omega_i)$.

\subsection{Convergence of the adaptive procedure}

In this section, we analyze the convergence of the above online enrichment procedure.
We begin our analysis at the inequality (\ref{eq:su_error1}). Notice that, this inequality can be written as
\begin{equation}
\| u_h - u_H^{(m+1)} \|_A^2 \leq  \| u_h - u_H^{(m)}\|_A^2 - \| R_i^{(m)} \|^2,
\label{eq:su_error2}
\end{equation}
%where the norm of the residual $\| R_i^{(m)} \|$ is defined as
%\begin{equation*}
%\| R_i^{(m)} \| = \sup_{v\in V_0^h(\omega_i)} \frac{| R_i^{(m)}(v)|}{\|  v\|_{\omega_i}},
%\end{equation*}
%and $\|v\|_{\omega_i}^2 = a_{\omega_i}(v,v)$.
when the basis function $\phi$ is obtained as in (\ref{eq:online}).
%by maximizing $R_i^{(m)}(\phi)$.

%Next, we define the following elliptic projection operator $P_{\text{DG}}: H^1(\mathcal{T}^H) \rightarrow V^h_{\text{DG}}$ by
%\begin{equation*}
%a_{\text{DG}}(P_{\text{DG}}w,v) = a_{\text{DG}}(w,v), \quad\forall v\in V^h_{\text{DG}}
%\end{equation*}
%where $H^1(\mathcal{T}^H)$ is the space of piecewise $H^1$ functions with respect to the coarse grid $\mathcal{T}^H$.
%It is then easy to show that $\|P_{\text{DG}} w\|_A \leq \|w\|_A$ for all $w\in H^1(\mathcal{T}^H)$.
%Notice the the projection operator preserves piecewise constant functions on the coarse grid $\mathcal{T}^H$. Consequently,
%\begin{equation}
%\sum_{K\in\mathcal{T}^H}\sum_{j=1}^4 P_{\text{DG}} \chi^{K}_j = 1.
%\label{eq:pou}
%\end{equation}
%Property (\ref{eq:pou}) holds since the sum $\sum_{j=1}^4 \chi^K_j = 1$ on each $K$,
%for all $K\in\mathcal{T}^H$.

On the other hand, we will show that the error $\| u_h - u_H^{(m)}\|_A$
can be controlled by the residual norm $\| R_i^{(m)}\|$. To do so,
we consider an arbitrary function $v\in V^h_{\text{DG}}$.
Let $v_i \in V^h(\omega_i)$ be the restriction of $v$ in $\omega_i$,
and let $v_i^{(0)}\in V^{i,\text{off}}$ be the component of $v_i$ in the offline space $V^{i,\text{off}}$.
By the GMsDGM (\ref{eq:ipdg}), the fine-grid problem (\ref{eq:ipdgfine}) and the
Galerkin orthogonality (\ref{eq:ortho}), we have
\begin{equation*}
a_{\text{DG}}(u_h - u_H^{(m)}, v) =  a_{\text{DG}}(u_h - u_H^{(m)}, v-v^{(0)}), \quad\forall v^{(0)} \in V_H^{(0)},
\end{equation*}
where we define $v^{(0)} = \sum_{i=1}^{N_c} v_i^{(0)} \in V_H^{(0)}$ and
use the fact that $V_H^{(0)} \subset V_H^{(m)}$ for all $m \geq 0$.
%Note that we can write $v^{(0)} = \sum_{i=1}^{N_c} v_i^{(0)}$, where $v_i^{(0)}$ is the component of $v^{(0)}$
%lying in $V^{i,\text{off}}$.
By (\ref{eq:ipdgfine}), we have
\begin{equation*}
a_{\text{DG}}(u_h - u_H^{(m)}, v) =  (f, v-v^{(0)}) - a_{\text{DG}}(u_H^{(m)}, v-v^{(0)}).
\end{equation*}
Using the property $\sum_{j=1}^4 \chi_{(j)}^K = 1$ for all $K\in\mathcal{T}^H$,
\begin{equation*}
a_{\text{DG}}(u_h - u^{(m)}_H, v) = \sum_{K\in\mathcal{T}^H} \sum_{j=1}^4  \Big( (f,\chi^K_{(j)} (v-v_i^{(0)})) -  a_{\text{DG}}(u^{(m)}_H, \chi^K_{(j)} (v-v_i^{(0)})) \Big).
\end{equation*}
Writing the above sum over coarse neighborhoods $\omega_i$, we have
\begin{equation*}
a_{\text{DG}}(u_h - u^{(m)}_H, v) = \sum_{i=1}^{N_c} \sum_{K\subset\omega_i} \Big( (f,\chi^K_{i} (v-v_i^{(0)})) -  a_{\text{DG}}(u^{(m)}_H, \chi^K_{i} (v-v_i^{(0)})) \Big).
\end{equation*}

For each coarse neighborhood $\omega_i$, we define the following modified local residual by
\begin{equation}
\widetilde{R}^{(m)}_i(v) = \sum_{K\subset\omega_i} \Big( (f,\chi_i^K v) - a_{\text{DG}}(u_H^{(m)},\chi_i^K v) \Big), \quad v \in V^h(\omega_i).
\label{eq:mres}
\end{equation}
The modified residual $\widetilde{R}_i^{(m)}$ can be seen as a linear functional defined on $V^h(\omega_i)$
with norm $\| \widetilde{R}_i^{(m)} \|$ defined in the following way
\begin{equation*}
\| \widetilde{R}_i^{(m)} \| = \sup_{v\in V^h(\omega_i)} \frac{| \widetilde{R}_i^{(m)}(v)|}{\|  \sum_{K\subset\omega_i} \chi_i^K v\|_{\omega_i}}
\end{equation*}
%and $\|v\|_{\omega_i}^2 = a_{\omega_i}(v,v)$.
In the above definitions, $\chi_i^K$ is considered to be defined only on $K$, and has zero value outside $K$.

Using the definition of the modified residual $\widetilde{R}_i^{(m)}$, we have
\begin{equation}
a_{\text{DG}}(u_h - u^{(m)}_H, v) \leq \sum_{i=1}^{N_c}  \| \widetilde{R}_i^{(m)}\| \, \|  \sum_{K\subset\omega_i} \chi_i^K (v-v_i^{(0)} ) \|_A
\label{eq:residualbound}
\end{equation}
where we used the fact that $\sum_{K\subset\omega_i} \chi_i^K (v-v_i^{(0)} )$ is zero on $\partial\omega_i$.
Using Lemma \ref{lem:norm},
\begin{equation}
\| \sum_{K\subset\omega_i} \chi_i^K (v-v_i^{(0)} ) \|_A \leq a_1^{\frac{1}{2}} \| \sum_{K\subset\omega_i} \chi_i^K (v-v_i^{(0)}) \|_{\text{DG}}.
\label{eq:norm}
\end{equation}
By the definition of the DG-norm,
\begin{equation}
\| \chi_i^K (v-v_i^{(0)}) \|_{\text{DG}}^2 = \sum_{K\subset \omega_i} \int_K \kappa |\nabla(\chi_i^K (v-v_i^{(0)}))|^2
+ \frac{\gamma}{h} \sum_{e} \int_e \overline{\kappa} \jump{\chi_i^K (v-v_i^{(0)})}^2.
\label{eq:b0}
\end{equation}
For each $K\subset\omega_i$, we have
\begin{equation}
\int_K \kappa |\nabla(\chi_i^K (v-v_i^{(0)}))|^2
\leq 2 \int_K \kappa \chi_i^2 |\nabla (v-v_i^{(0)})|^2 + 2\int_K \kappa |\nabla \chi_i^K|^2 (v-v_i^{(0)})^2.
\label{eq:b1}
\end{equation}
For each $e \in\mathcal{E}^H_i$, we have
\begin{equation}
\int_e \overline{\kappa} \jump{\chi_i^K (v-v_i^{(0)})}^2
\leq 2 \int_e \overline{\kappa} \average{\chi_i^K}^2 \jump{v-v_i^{(0)}}^2 + 2 \int_e \overline{\kappa} \jump{ \chi_i^K}^2 \average{v-v_i^{(0)}}^2.
\label{eq:b2}
\end{equation}
Combining inequalities (\ref{eq:b1}) and (\ref{eq:b2}) in (\ref{eq:b0}), we have
\begin{equation}
\| \chi_i^K (v-v_i^{(0)}) \|_{\text{DG}}^2 \leq 2 \| v-v_i^{(0)}\|_{A_i}^2
+ 2 \Big( \int_{\omega_i} \overline{\kappa} |\nabla \chi_i^K |^2 (v-v_i^{(0)})^2 + \frac{\gamma}{h} \sum_{e\in\mathcal{E}^H_i} \int_e \overline{\kappa} \jump{ \chi_i^K}^2 \average{v-v_i^{(0)}}^2 \Big)
\label{eq:b3}
\end{equation}
where $\|v\|_{A_i}^2 = a_{\omega_i}(v,v)$.
Using the spectral problem (\ref{eq:spec}), we have
\begin{equation*}
\| v-v_i^{(0)}\|_{A_i}^2 = a_{\omega_i}(v_i-v_i^{(0)}, v_i-v_i^{(0)}) \leq a_{\omega_i}(v_i, v_i) = \|v\|_{A_i}^2
\end{equation*}
and
\begin{equation*}
\int_{\omega_i} \overline{\kappa} |\nabla \chi_i^K |^2 (v-v_i^{(0)})^2 + \frac{\gamma}{h} \sum_{e\in\mathcal{E}^H_i} \int_e \overline{\kappa} \jump{\nabla \chi_i^K}^2 \average{v-v_i^{(0)}}^2
= s_{\omega_i} (v_i-v_i^{(0)}, v_i-v_i^{(0)})
\leq \frac{1}{\lambda_{L_i+1}^{\omega_i}} \| v\|_{A_i}^2.
\end{equation*}
Thus, (\ref{eq:b3}) and (\ref{eq:norm}) impleis
\begin{equation*}
\| \chi_i^K (v-v_i^{(0)}) \|_{A}^2 \leq 2 a_1 \Big( 1 + \frac{1}{\lambda^{\omega_i}_{L_i+1}} \Big)  \| v\|_{A_i}^2.
\end{equation*}
%Thus,
%\begin{equation*}
%\| \chi_j^K (v-v_i^{(0)} ) \|_A^2 \leq 2 a_1 ( 1 + \frac{1}{\lambda^{\omega_i}_k} )  \| v\|_{A_i}^2.
%\end{equation*}

Hence, (\ref{eq:residualbound}) becomes
\begin{equation*}
a_{\text{DG}}(u_h-u_H^{(m)},v) \leq \Big( \sum_{i=1}^{N_c} 2 a_1 ( 1 + \frac{1}{\lambda^{\omega_i}_{L_i+1}} ) \|\widetilde{R}_i^{(m)}\|^2 \Big)^{\frac{1}{2}} \, \Big(\sum_{i=1}^{N_c} \| v\|_{A_i}^2\Big)^{\frac{1}{2}}.
\end{equation*}
We remark that the above inequality holds for any $v\in V^h_{\text{DG}}$. Taking $v=v_h-v_H^{(m)}$ and using Lemma \ref{lem:norm},
we finally obtain
\begin{equation}
\| u_h - u_H^{(m)} \|_A^2 \leq 2a_0^{-1} a_1 C_0 \sum_{i=1}^{N_c}  \Big( 1 + \frac{1}{\lambda^{\omega_i}_{L_i+1}} \Big) \|\widetilde{R}_i^{(m)}\|^2,
\label{eq:residualbound1}
\end{equation}
where $C_0 = \max_{K\in\mathcal{T}^H} n_K$ and $n_K$ is the number of vertices of the coarse grid block $K$.

We define
\begin{equation}
\theta = \|R_i^{(m)}\|^2 / \eta^2, \quad\text{and}\quad
\eta^2 = 2a_0^{-1} a_1 C_0  \sum_{i=1}^{N_c}  \Big( 1 + \frac{1}{\lambda^{\omega_i}_{L_i+1}} \Big) \|\widetilde{R}_i^{(m)}\|^2.
\label{eq:theta}
\end{equation}
From (\ref{eq:su_error2}) and (\ref{eq:residualbound1}), we see that
the following convergence holds
\begin{equation*}
\| u_h - u_H^{(m+1)}\|_A^2 \leq ( 1 - \theta ) \| u_h - u_H^{(m)}\|_A^2.
\end{equation*}

We summarize the above results in the following theorem.

\begin{theorem}
Let $u_h$ be the solution of (\ref{eq:ipdgfine}) and $u_H^{(m)}$, $m\geq 0$, be the solution of (\ref{eq:ipdg}) in the $m$-th iteration. Then the following residual bound holds
\begin{equation}
\| u_h - u_H^{(m)} \|_A^2 \leq 2a_0^{-1} a_1 C_0 \sum_{i=1}^{N_c}  \Big( 1 + \frac{1}{\lambda^{\omega_i}_{L_i+1}} \Big) \|\widetilde{R}_i^{(m)}\|^2.
\label{eq:apost} 
\end{equation}
Moreover, the following convergence holds
\begin{equation}
\| u_h - u_H^{(m+1)}\|_A^2 \leq ( 1 -\theta ) \| u_h - u_H^{(m)}\|_A^2
\label{eq:adaptive_conv} 
\end{equation}
where $\theta$ is defined in (\ref{eq:theta}).
\end{theorem}

We remark that one can derive a priori error estimate for the error $\|u_h-u_H^{(m)}\|_{\text{DG}}$, for every $m\geq 0$.
Since the purpose of this paper is an a posteriori error estimate (\ref{eq:apost}) and the convergence of an adaptive enrichment algorithm (\ref{eq:adaptive_conv}),
we will not derive a priori error estimate.

Finally, we remark that by using more basis functions in the initial space $V_H^{(0)}$,
the values of the eigenvalues $\lambda_{L_i+1}^{\omega_i}$ are larger.
Thus, the value of $\theta$ is further away from zero, and this fact enhances the convergence rate.
In particular, the convergence rate is affected by the quantity $\Lambda_{\text{min}} = \min_{1\leq i\leq N_c} \lambda_{L_i+1}^{\omega_i}$.
The convergence is slow when $\Lambda_{\text{min}}$ is small (cf. \cite{cel15_1, ge09_1reduceddim}).
We also remark that one can add online basis functions in multiple coarse neighborhoods to speed up the convergence.
Let $S$ be the index set for which online basis functions are added in $\omega_i$ for $i\in S$.
By using similar arguments as above, we obtain
\begin{equation*}
\| u_h - u_H^{(m+1)}\|_A^2 \leq ( 1 -\widetilde{\theta} ) \| u_h - u_H^{(m)}\|_A^2
\end{equation*}
where
\begin{equation*}
\widetilde{\theta} = \sum_{i\in S} \| R_i^{(m)}\|^2 / \eta^2.
\end{equation*}

% Numerical
\section{Numerical Results}
\label{sec:numerresults}

In this section, we will present some numerical examples to
show the performance of the
proposed method.
The implementation procedure of online adaptive GMsDGM is described below.
First, we choose a fixed number of functions
for every coarse neighborhood by solving the local spectral problem. This fixed number for every coarse neighborhood is called the number of initial basis. After that, we split these functions into the basis functions of the offline space such that each basis function is supported in one coarse grid block. We denote this offline space as $V^{\text{off}}$ and set $V_H^{(0)} = V^{\text{off}}$.

The coarse neighborhoods are denoted by $\omega_{i,j}$,
where $i=1,2,\cdots, N_x$ and $j=1,2,\cdots, N_y$
and $N_x$ and $N_y$ are the number of coarse nodes in the $x$ and $y$ directions respectively.
We consider $I_{x,\text{odd}}$ and $I_{x,\text{even}}$ as the set of odd and even indices from $\{ 1,2,\cdots, N_x\}$.
Similarly, $I_{y,\text{odd}}$ and $I_{y,\text{even}}$ are the set of odd and even indices from $\{ 1,2,\cdots, N_y\}$.
In each iteration of our online adaptive GMsDGM, we will perform $4$ sub-iterations which add online basis functions
in the non-overlapping coarse neighborhoods $\omega_{i,j}$ with $(i,j) \in I_{x,\text{odd}} \times I_{y,\text{odd}}$,
$(i,j) \in I_{x,\text{odd}} \times I_{y,\text{even}}$, $(i,j) \in I_{x,\text{even}} \times I_{y,\text{odd}}$
and $(i,j) \in I_{x,\text{even}} \times I_{y,\text{even}}$ respectively.

\begin{figure}[htb]
  \centering
  \includegraphics[width=0.40 \textwidth]{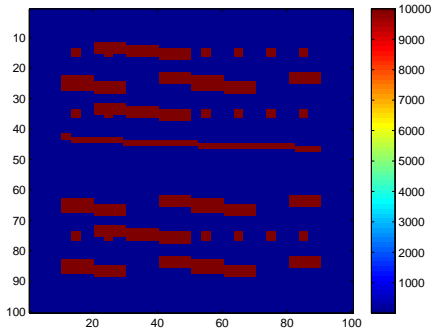}
  \caption{Permeability field $\kappa$.}
  \label{fig:medium and source}
\end{figure}

We will take $\gamma=2$ and $D=[0,1]^2$. The domain is divided into
$10\times 10$ uniform square coarse blocks. Each coarse block is then divided
into $10\times 10$ fine blocks consisting of uniform squares. Namely, the whole domain
is partitioned by $100\times 100$ fine grid blocks.
The medium parameter $\kappa$ is shown in Figure
\ref{fig:medium and source}. The source function $f$ is taken as the constant 1.
To compare the accuracy, we will use the following error
quantities
\[
e_{2}=\cfrac{\|u_h - u_H \|_{L^{2}(D)}}{\|u_h\|_{L^{2}(D)}},\quad\text{and}\quad  e_{a}=\cfrac{\| u_h - u_H\|_{\text{DG}}}{\| u_h \|_{\text{DG}}}.
\]

\subsection{Comparison of using different number of initial basis}

In Table \ref{tab:fullspace_1 2_3_4basis}, we present the convergence history of
our algorithm for using one, two, three, four initial basis per coarse neighborhood.
Notice that, in
the presentation of our results, DOF means the total number of basis functions used in the whole domain.
We use the continuous multiscale basis functions as the initial
partition of unity.
In the tables, we obtain a fast error decay which give us a numerical solution with error smaller than $0.1\%$ in two or three iterations. We can see the error decay of using one initial basis is slower than the error decay of using two or more initial basis since $\Lambda_{\text{min}}$ for using one initial basis is too small.

\begin{table}[ht]
\begin{centering}
\begin{tabular}{|c|c|c|}
\hline
DOF & $e_{a}$ & $e_{2}$\tabularnewline
\hline
324 & 44.50\% & 24.88\%\tabularnewline
\hline
648 & 9.92\% & 2.18\%\tabularnewline
\hline
972 & 0.78\% & 7.54e-2\%\tabularnewline
\hline
1296 & 3.24e-2\% & 2.13e-3\%\tabularnewline
\hline
1620 & 2.42e-4\% & 1.10e-5\%\tabularnewline
\hline
\end{tabular}
\begin{tabular}{|c|c|c|}
\hline
DOF & $e_{a}$ & $e_{2}$\tabularnewline
\hline
648 & 17.73\% & 3.58\%\tabularnewline
\hline
972 & 0.31\% & 1.80e-2\%\tabularnewline
\hline
1296 & 3.52e-3\% & 1.62e-4\%\tabularnewline
\hline
1620 & 1.81e-5\% & 8.58e-7\%\tabularnewline
\hline
1948 & 1.04e-7\% & 4.68e-9\%\tabularnewline
\hline
\end{tabular}

\begin{tabular}{|c|c|c|}
\hline
DOF & $e_{a}$ & $e_{2}$\tabularnewline
\hline
972 & 11.30\% & 1.72\%\tabularnewline
\hline
1296 & 0.45\% & 2.44e-2\%\tabularnewline
\hline
1620 & 3.05e-3\% & 1.37e-4\%\tabularnewline
\hline
1944 & 1.06e-5\% & 4.08e-7\%\tabularnewline
\hline
2240 & 4.59e-8\% & 2.14e-9\%\tabularnewline
\hline
\end{tabular}
\begin{tabular}{|c|c|c|}
\hline
DOF & $e_{a}$ & $e_{2}$\tabularnewline
\hline
1296 & 8.38\% & 1.00\%\tabularnewline
\hline
1620 & 7.98e-2\% & 3.13e-3\%\tabularnewline
\hline
1944 & 9.93e-4\% & 3.57e-5\%\tabularnewline
\hline
2268 & 1.39e-5\% & 5.15e-7\%\tabularnewline
\hline
2540 & 4.23e-8\% & 1.55e-9\%\tabularnewline
\hline
\end{tabular}
\par\end{centering}
%\centering{}
\protect\caption{Top-left: One initial basis ($\Lambda_{\text{min}}=4.89e-4$). Top-right: Two initial basis ($\Lambda_{\text{min}}=0.9504$). \protect \\ Bottom-left: Three initial basis ($\Lambda_{\text{min}}=1.4226$). Bottom-right: Four initial basis ($\Lambda_{\text{min}}=2.2045$).}
\label{tab:fullspace_1 2_3_4basis}
\end{table}

\begin{figure}[ht]
  \centering
  \includegraphics[width=0.40 \textwidth]{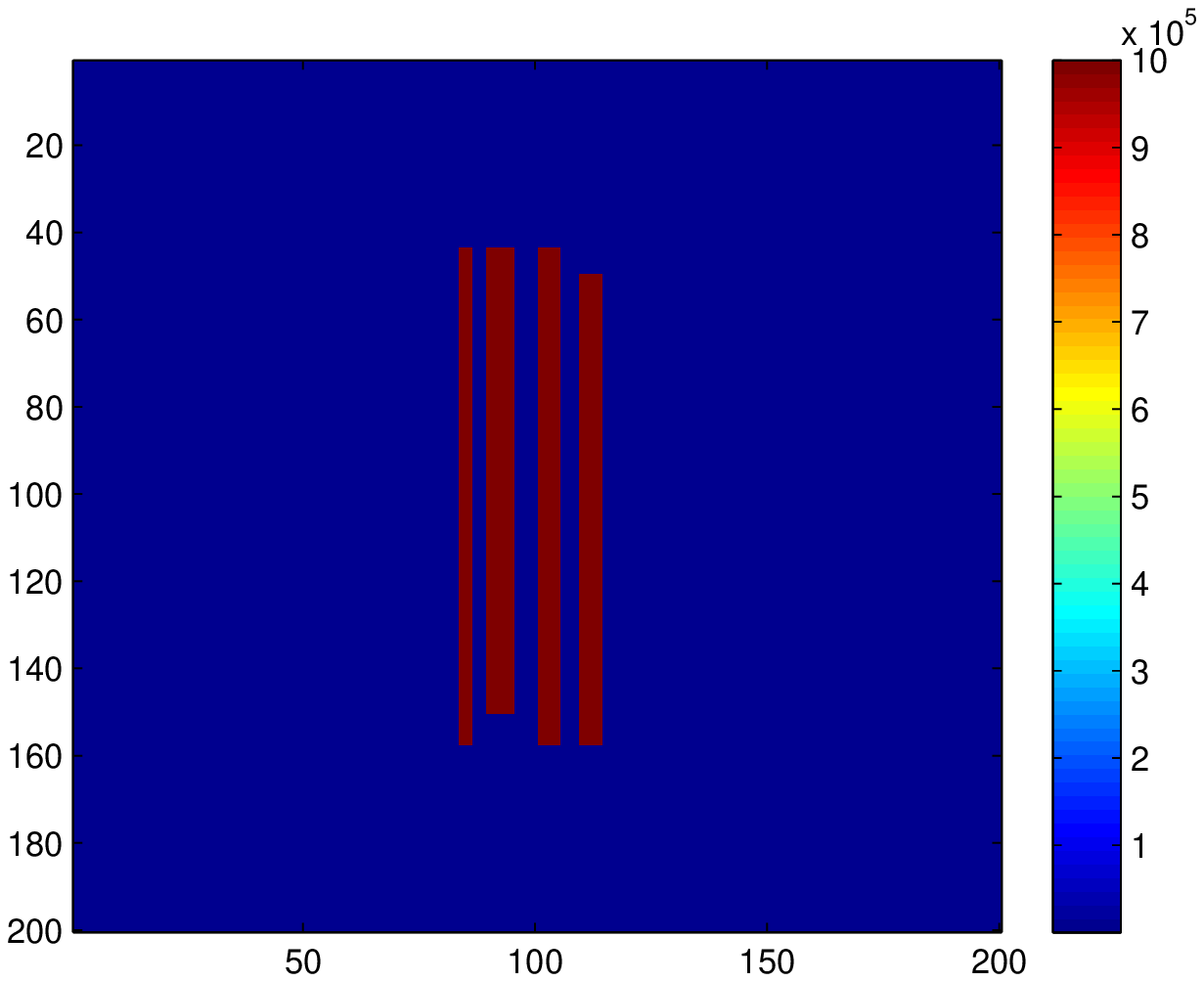}
  \caption{Permeability field $\kappa$.}
  \label{fig:medium2}
\end{figure}

%\begin{figure}[htb]
%  \centering
%  \includegraphics[width=0.4 \textwidth]{error_compare_tilde}
%  \includegraphics[width=0.4 \textwidth]{error_compare_1e6}
%  \caption{Error comparison.  Along $x$-axis: Dimensions of $V_{\text{ms}}$. Along $y$-axis: Relative energy errors. Left: $1e4$. Right: $1e6$.}
%  \label{fig:compare error}
%\end{figure}
To further study the importance of the initial basis, we will present another example with a different medium parameter $\kappa$ shown in Figure \ref{fig:medium2}. The domain $D$ is divided into $5\times 5$ coarse blocks consisting of uniform squares. Each coarse block is then divided into $40\times 40$ fine blocks also consisting of uniform squares. The convergence history for the use of one, two, three, four initial basis and the corresponding total number of degrees of freedom (DOF)
are shown in Table
\ref{fullspace_onebasis_4channels}, Table
\ref{fullspace_twobasis_4channels}, Table
\ref{fullspace_threebasis_4channels}, Table
\ref{fullspace_fourbasis_4channels} respectively.
We consider two different contrasts. On the right table, we increase the contrast by 100 times. More precisely, the
conductivity of inclusions and channels in Figure 2 (left figure) is multiplied by 100.
In this case, the first 4 eigenvalue that are in the regions with channels become 100 times smaller. The decrease in the eigenvalues will slow down the error decay. In Table \ref{fullspace_onebasis_4channels}, we can observe that the error decay for the lower contrast case is much faster than the higher contrast case. In the higher contrast case, the error stop decreasing in some iterations. Similar observations are obtained when we use 2 or 3 initial basis. For using four initial basis, we observe a rapid convergence for both higher and lower contrast case.

%\begin{figure}[htb]
%  \centering
%  \includegraphics[width=0.65 \textwidth]{channel_error_tilde}
%  \caption{Error comparison for different number of initial basis functions. Along $x$-axis: Dimensions of $V_{\text{ms}}$. Along $y$-axis: Relative energy errors.}
%  \label{fig:compare error_channels}
%\end{figure}

%In conclusion, we observe
%\begin{itemize}
%
%\item If $V_{\text{ms}}$ does not satisfy ONERP, then the error
%decay is slower as the contrast becomes larger.
%
%\item If $V_{\text{ms}}$ does not satisfy ONERP,
%in some cases, we have observed the error does not decrease
%as we add online basis functions (see Table \ref{fullspace_onebasis_4channels}, \ref{tab:fullspace_2_3_4basis_4channels}).
%
%
%\item If $V_{\text{ms}}$ satisfies ONERP, then we observe a fast
%convergence, which is independent of contrast.
%
%
%\end{itemize}

%{\bf Our conclusion of Section 4.1.1. (1) We see from Table 1 that if we choose one basis, then the reduction is small. This is because of $\lambda$ term. Thus, one needs to start with several basis functions with $\Lambda$ being larger than 0. (2) The error decay is much faster with online basis vs. offline after several basis. (3) As we see that if we choose 2 basis, with one iteration we can get an error less than 1 \% while for 3 basis, even though there is larger error reduction, the extra basis function in the offline space is not necessary. Indeed, from permeability, we know that $\Lambda>0$ if we have at least two basis functions because there are few blocks with two channels.}

\begin{table}[ht]
\begin{centering}
\begin{tabular}{|c|c|c|}
\hline
DOF & $e_{a}$ & $e_{2}$\tabularnewline
\hline
64 & 25.44\% & 6.67\%\tabularnewline
\hline
128 & 1.20\% & 0.23\%\tabularnewline
\hline
192 & 0.47\% & 0.10\%\tabularnewline
\hline
256 & 0.26\% & 5.79e-2\%\tabularnewline
\hline
320 & 0.10\% & 2.30e-2\%\tabularnewline
\hline
384 & 6.22e-2\% & 1.02e-2\%\tabularnewline
\hline
448 & 3.70e-4\% & 1.57e-5\%\tabularnewline
\hline
\end{tabular} %
\begin{tabular}{|c|c|c|}
\hline
DOF & $e_{a}$ & $e_{2}$\tabularnewline
\hline
64 & 25.45\% & 6.67\%\tabularnewline
\hline
128 & 1.45\% & 0.27\%\tabularnewline
\hline
192 & 1.39\% & 0.27\%\tabularnewline
\hline
256 & 0.84\% & 0.15\%\tabularnewline
\hline
320 & 0.34\% & 7.98e-2\%\tabularnewline
\hline
384 & 0.34\% & 7.91e-2\%\tabularnewline
\hline
448 & 0.15\% & 3.71e-2\%\tabularnewline
\hline
\end{tabular}
\par\end{centering}

\protect\caption{One initial basis. Left: Lower contrast(1e4)($\Lambda_{\text{min}}=0.0062$).  \protect\\
Right: Higher contrast(1e6)($\Lambda_{\text{min}}=6.22e-5$).}
\label{fullspace_onebasis_4channels}
\end{table}

\begin{table}[ht]
\begin{centering}
\begin{tabular}{|c|c|c|}
\hline
DOF & $e_{a}$ & $e_{2}$\tabularnewline
\hline
128 & 18.22\% & 4.42\%\tabularnewline
\hline
192 & 1.14\% & 0.12\%\tabularnewline
\hline
256 & 0.50\% & 4.95e-2\%\tabularnewline
\hline
320 & 4.17e-2\% & 2.06e-3\%\tabularnewline
\hline
384 & 5.73e-3\% & 5.38e-4\%\tabularnewline
\hline
448 & 7.12e-4\% & 2.89e-5\%\tabularnewline
\hline
\end{tabular} %
\begin{tabular}{|c|c|c|}
\hline
DOF & $e_{a}$ & $e_{2}$\tabularnewline
\hline
128 & 18.56\% & 4.62\%\tabularnewline
\hline
192 & 1.37\% & 0.16\%\tabularnewline
\hline
256 & 1.25\% & 0.14\%\tabularnewline
\hline
320 & 1.23\% & 0.13\%\tabularnewline
\hline
384 & 0.41\% & 3.22e-2\%\tabularnewline
\hline
448 & 3.63e-2\% & 3.56e-3\%\tabularnewline
\hline
\end{tabular}
\par\end{centering}
\protect\caption{Two initial basis. Left: Lower contrast(1e4)($\Lambda_{\text{min}}=0.027$).  \protect\\
Right: Higher contrast(1e6)($\Lambda_{\text{min}}=2.72e-4$).}
\label{fullspace_twobasis_4channels}
\end{table}

\begin{table}[ht]
\begin{centering}
\begin{tabular}{|c|c|c|}
\hline
DOF & $e_{a}$ & $e_{2}$\tabularnewline
\hline
192 & 10.69\% & 1.86\%\tabularnewline
\hline
256 & 0.80\% & 6.66e-2\%\tabularnewline
\hline
320 & 0.34\% & 2.24e-2\%\tabularnewline
\hline
384 & 1.51e-2\% & 6.24e-4\%\tabularnewline
\hline
448 & 2.25e-4\% & 1.61e-5\%\tabularnewline
\hline
508 & 1.72e-6\% & 6.70e-8\%\tabularnewline
\hline
\end{tabular} %
\begin{tabular}{|c|c|c|}
\hline
DOF & $e_{a}$ & $e_{2}$\tabularnewline
\hline
192 & 11.55\% & 2.14\%\tabularnewline
\hline
256 & 1.13\% & 0.10\%\tabularnewline
\hline
320 & 0.98\% & 8.85e-2\%\tabularnewline
\hline
384 & 0.96\% & 8.95e-2\%\tabularnewline
\hline
448 & 0.30\% & 1.39e-2\%\tabularnewline
\hline
508 & 2.00e-3\% & 8.39e-5\%\tabularnewline
\hline
\end{tabular}
\par\end{centering}
\protect\caption{Three initial basis. Left: Lower contrast(1e4)($\Lambda_{\text{min}}=0.0371$). \protect\\
Right: Higher contrast(1e6)($\Lambda_{\text{min}}=3.75e-4$).}
\label{fullspace_threebasis_4channels}
\end{table}

\begin{table}[ht]
\begin{centering}
\begin{tabular}{|c|c|c|}
\hline
DOF & $e_{a}$ & $e_{2}$\tabularnewline
\hline
248 & 7.92\% & 1.14\%\tabularnewline
\hline
312 & 0.25\% & 2.42e-2\%\tabularnewline
\hline
376 & 5.09e-3\% & 2.72e-4\%\tabularnewline
\hline
440 & 5.18e-5\% & 2.62e-6\%\tabularnewline
\hline
484 & 1.39e-6\% & 6.40e-8\%\tabularnewline
\hline
\end{tabular} %
\begin{tabular}{|c|c|c|}
\hline
DOF & $e_{a}$ & $e_{2}$\tabularnewline
\hline
242 & 9.63\% & 1.59\%\tabularnewline
\hline
306 & 0.51\% & 5.40e-2\%\tabularnewline
\hline
370 & 1.38e-2\% & 9.46e-4\%\tabularnewline
\hline
434 & 2.10e-4\% & 1.59e-5\%\tabularnewline
\hline
494 & 1.74e-6\% & 1.27e-7\%\tabularnewline
\hline
\end{tabular}
\par\end{centering}
\protect\caption{Four initial basis. Left: Lower contrast(1e4)($\Lambda_{\text{min}}=0.4472$).  \protect\\
Right: Higher contrast(1e6)($\Lambda_{\text{min}}=0.3844$).}
\label{fullspace_fourbasis_4channels}
\end{table}

%\clearpage{}

\subsection{Setting tolerance for the residual}

In this section, we will show the performance for the online enrichment implementing it only for regions
with a residual error bigger than a certain threshold. We consider the medium parameter shown in Figure \ref{fig:medium and source}. We show the results for using three different tolerances ($tol$)
$10^{-3}$, $10^{-4}$ and $10^{-5}$. We will enrich for the coarse regions with
residual larger than the tolerance.
In Table \ref{tab:tol one basis},
we show the errors when using $1$ initial basis function
for tolerances $10^{-3}$, $10^{-4}$ and $10^{-5}$.
We can see that the convergence history in the first few iteration is similar to the result shown in previous section. Moreover, the energy error of the multiscale solution
is in the same order of the tolerance and the error will stop decreasing even if we perform more iterations. Therefore, we can compute a multiscale solution with a prescribed error level
by choosing a suitable tolerance in the adaptive algorithm.
In Table  \ref{tab:tol two basis} and Table \ref{tab:tol three basis},
we show the errors for the last three iterations
when using $2$ and $3$ initial basis functions respectively
for tolerances $10^{-3}$, $10^{-4}$ and $10^{-5}$.
We have the same observation that the energy errors have the same magnitude as the tolerances.

%In the second case,
%the online enrichment is performed for
%coarse regions that have cumulative residual that is $\theta$ fraction
%of the total residual. One of our objectives is to show that one can drive
%the error down to a number below a threshold, adaptively.

%Next, we will present some numerical examples showing controlling the local residual error can control the finial energy error. More precisely, we set a tolerance such that we will not add the residual based local basis into the finite element space if the relative energy norm of that basis is smaller than the tolerance and we will stop adding basis if all of the relative energy norm of the basis are smaller than the tolerance.

\begin{table}[htb]
\begin{centering}
\begin{tabular}{|c|c|c|}
\hline
DOF & $e_{a}$ & $e_{2}$\tabularnewline
\hline
324 & 44.50\% & 24.88\%\tabularnewline
\hline
648 & 9.92\% & 2.18\%\tabularnewline
\hline
924 & 0.81\% & 7.72e-2\%\tabularnewline
\hline
976 & 0.29\% & 2.49e-2\%\tabularnewline
\hline
\end{tabular}
\begin{tabular}{|c|c|c|}
\hline
DOF & $e_{a}$ & $e_{2}$\tabularnewline
\hline
324 & 44.50\% & 24.88\%\tabularnewline
\hline
648 & 9.92\% & 2.18\%\tabularnewline
\hline
972 & 0.78\% & 7.54e-2\%\tabularnewline
\hline
1176 & 4.12e-2\% & 2.88e-3\%\tabularnewline
\hline
1184 & 2.65e-2\% & 1.57e-3\%\tabularnewline
\hline
\end{tabular}
\begin{tabular}{|c|c|c|}
\hline
DOF & $e_{a}$ & $e_{2}$\tabularnewline
\hline
324 & 44.50\% & 24.88\%\tabularnewline
\hline
648 & 9.92\% & 2.18\%\tabularnewline
\hline
972 & 0.78\% & 7.54e-2\%\tabularnewline
\hline
1284 & 3.24e-2\% & 2.13e-3\%\tabularnewline
\hline
1364 & 2.56e-3\% & 1.55e-4\%\tabularnewline
\hline
\end{tabular}
\par\end{centering}

\centering{}\protect\caption{One initial basis. Left: $tol = 10^{-3}$. Middle: $tol = 10^{-4}$. Right: $tol = 10^{-5}$.}
\label{tab:tol one basis}
\end{table}

\begin{table}[htb]
\begin{centering}
\begin{tabular}{|c|c|c|}
\hline
DOF & $e_{a}$ & $e_{2}$\tabularnewline
\hline
648 & 17.73\% & 3.58\%\tabularnewline
\hline
964 & 0.33\% & 1.85e-2\%\tabularnewline
\hline
972 & 0.30\% & 1.63e-2\%\tabularnewline
\hline
\end{tabular} %
\begin{tabular}{|c|c|c|}
\hline
DOF & $e_{a}$ & $e_{2}$\tabularnewline
\hline
648 & 17.73\% & 3.58\%\tabularnewline
\hline
972 & 0.31\% & 1.80e-2\%\tabularnewline
\hline
1136 & 2.53e-2\% & 1.24e-3\%\tabularnewline
\hline
\end{tabular} %
\begin{tabular}{|c|c|c|}
\hline
DOF & $e_{a}$ & $e_{2}$\tabularnewline
\hline
972 & 0.31\% & 1.80e-2\%\tabularnewline
\hline
1248 & 3.99e-3\% & 1.85e-4\%\tabularnewline
\hline
1276 & 2.49e-3\% & 1.19e-4\%\tabularnewline
\hline
\end{tabular}
\par\end{centering}

\centering{}\protect\caption{Two initial basis. Left: $tol = 10^{-3}$. Middle: $tol = 10^{-4}$. Right: $tol = 10^{-5}$.}
\label{tab:tol two basis}
\end{table}

\begin{table}[htb]
\begin{centering}
\begin{tabular}{|c|c|c|}
\hline
DOF & $e_{a}$ & $e_{2}$\tabularnewline
\hline
972 & 11.30\% & 1.72\%\tabularnewline
\hline
1248 & 0.50\% & 2.57e-2\%\tabularnewline
\hline
1276 & 0.24\% & 9.98e-3\%\tabularnewline
\hline
\end{tabular} %
\begin{tabular}{|c|c|c|}
\hline
DOF & $e_{a}$ & $e_{2}$\tabularnewline
\hline
972 & 11.30\% & 1.72\%\tabularnewline
\hline
1296 & 0.45\% & 2.44e-2\%\tabularnewline
\hline
1436 & 2.60e-2\% & 9.70e-4\%\tabularnewline
\hline
\end{tabular} %
\begin{tabular}{|c|c|c|}
\hline
DOF & $e_{a}$ & $e_{2}$\tabularnewline
\hline
1296 & 0.45\% & 2.44e-2\%\tabularnewline
\hline
1564 & 3.52e-3\% & 1.56e-4\%\tabularnewline
\hline
1576 & 2.49e-3\% & 1.04e-4\%\tabularnewline
\hline
\end{tabular}
\par\end{centering}

\centering{}\protect\caption{Three initial basis. Left: $tol = 10^{-3}$. Middle: $tol = 10^{-4}$. Right: $tol = 10^{-5}$.}
\label{tab:tol three basis}
\end{table}

\subsection{Adaptive online enrichment}

In this section, we will show the performance for the online enrichment
implementing it only for regions
that have a cumulative residual that is $\theta$ fraction
of the total residual.
% and $H^{-1}$-norm of the residual larger than the tolerance.
We consider the medium parameter shown in Figure \ref{fig:medium2} (4 channels medium).

%More precisely, we consider
%$r_i = \|\phi_i\|_a$ where $\phi_i$ are the
%local basis function constructed by our algorithm.
Assume that the local residuals are arranged such that
\[
r_1\geq r_2 \geq r_3 \geq \cdots.
\]
We only add the basis {$\phi_1, \cdots, \phi_k$} for the coarse neighborhoods $\omega_1,\cdots, \omega_k$
such that $k$ is the smallest integer with
\[
\theta \sum^{N_c}_{i=1}r^2_i\leq \sum_{i=1}^k r^2_i.
\]

\begin{figure}[htb]
  \centering
  \includegraphics[width=0.40 \textwidth]{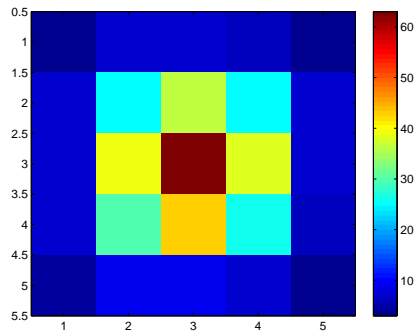}
  \caption{Distribution of number of basis functions in coarse blocks.}
  \label{fig:number of basis}
\end{figure}

In Table \ref{tab:ad tol1}, we present the error for the last $5$ iterations
when using $1$
initial basis functions with the tolerance $10^{-5}$ and $\theta=0.5$.
Comparing the result to the previous case, we can observe that this can use less number of basis functions to achieve a similar error.
In Figure~\ref{fig:number of basis}, we present the distribution of number of basis functions in coarse blocks,
and see that the number of basis functions is larger near the channels (c.f. Figure~\ref{fig:medium2}).
Thus,
online basis functions can be adaptively added in some regions using
an error indicator.

\begin{table}[htb]
\begin{centering}
\begin{tabular}{|c|c|c|}
\hline
DOF & $e_{a}$ & $e_{2}$\tabularnewline
\hline
348 & 0.35\% & 5.19e-2\%\tabularnewline
\hline
368 & 0.27\% & 4.03e-2\%\tabularnewline
\hline
392 & 6.13e-2\% & 9.34e-3\%\tabularnewline
\hline
412 & 6.04e-3\% & 6.60e-4\%\tabularnewline
\hline
424 & 1.51e-3\% & 1.25e-4\%\tabularnewline
\hline
\end{tabular}

\par\end{centering}

\centering{}\protect\caption{The results using cumulative errors with $\theta=0.5$, $tol=10^{-5}$ and 1 initial basis.}
\label{tab:ad tol1}
\end{table}

%In figure \ref{fig:number of basis}, we will show the dimension of local offline space in each coarse grid block.

\section{Conclusion}
\label{sec:conclusion}

Though the use of offline basis functions is important for multiscale finite
element methods, adding online basis functions in some regions can
improve the convergence dramatically. The construction of online basis functions
for various applications and discretizations require a careful analysis.
In particular, as we have shown earlier \cite{cel15_1} for GMsFEM within continuous
Galerkin framework that one needs a certain number of offline basis functions
in order to guarantee that the online basis functions can result to a convergence independent of physical parameters.
In this paper, we develop an online basis procedure for GMsDGM
that can provide a convergence independent of the contrast and small scales.
Because  multiscale basis functions are discontinuous across coarse-grid
boundaries, we construct a special offline space as well as
 online basis functions.
We show that our construction will guarantee a convergence independent
of the contrast and small scales if we select a certain number of
offline basis functions based on a local spectral problem.
Furthermore, we apply an adaptive procedure to add online basis functions
in only some selected regions. Numerical results are presented
to back up our theoretical findings.

\section*{Appendix}

In this section, we proof Lemma \ref{lem:norm}.
Let $K$ be a coarse grid block and let ${n}_{\partial K}$ be the
unit outward normal vector on $\partial K$.
We let $V^h(K)$ be the space of continuous piecewise bi-quadratic polynomials defined in $K$, and
we denote $V^{h}(\partial K)$
by the restriction of the conforming space $V^{h}(K)$ on $\partial K$.
The normal flux $\kappa \nabla u \cdot \,{n}_{\partial K}$ is understood as
an element in $V_{h}(\partial K)$ and is defined by
\begin{equation}
\int_{\partial K}({\kappa}\nabla{u}\cdot{n}_{\partial K})\cdot{v}=\int_{K} \kappa\nabla u\cdot\nabla\widehat{{v}},\quad{v}\in V^{h}(\partial K),\label{eq:flux}
\end{equation}
where $\widehat{{v}} \in V_h(K)$ is the harmonic extension of ${v}$ in $K$.
By the Cauchy-Schwarz inequality,
\[
\int_{\partial K}({\kappa}\nabla{u}\cdot{n}_{\partial K})\cdot{v}\leq a_{H}^{K}(u,u)^{\frac{1}{2}}\, a_{H}^{K}(\widehat{v},\widehat{v})^{\frac{1}{2}}.
\]
By an inverse inequality and the fact that $\widehat{v}$ is the harmonic
extension of $v$
\begin{equation}
a_{H}^{K}(\widehat{v},\widehat{v})\leq\kappa_{K}C_{\text{inv}}^{2}h^{-1}\int_{\partial K}|v|^{2},
\label{eq:maxeig}
\end{equation}
where we recall that $\kappa_{K}$ is the maximum of $\kappa$ over $K$ and $C_{\text{inv}}>0$ is
the constant from inverse inequality. Thus,
\[
\int_{\partial K}({\kappa}\nabla{u}\cdot{n}_{\partial K})\cdot{v}\leq\kappa_{K}^{\frac{1}{2}}C_{\text{inv}}h^{-\frac{1}{2}}\|v\|_{L^{2}(\partial K)}\, a_{H}^{K}(u,u)^{\frac{1}{2}}.
\]
This shows that
\begin{equation}
\int_{\partial K}|{\kappa}\nabla{u}\cdot{n}_{\partial K}|^{2}\leq\kappa_{K}C_{\text{inv}}^{2}h^{-1}a_{H}^{K}(u,u).
\label{eq:fluxbound}
\end{equation}

Next, by the definition of $a_{\text{DG}}$,
we have
\[
a_{\text{DG}}({u},{v})=a_{H}({u},{v})-\sum_{E\in\mathcal{E}^{H}}\int_{E}\Big(\average{{\kappa}\nabla{u}\cdot{n}_{E}}\jump{v}+\average{{\kappa}\nabla{v}\cdot{n}_{E}}\jump{u}\Big)+\sum_{E\in\mathcal{E}^{H}}\frac{\gamma}{h}\int_{E}\overline{\kappa}\jump{u} \jump{v}.
\]
Notice that
\[
a_{H}({u},{v})+\sum_{E\in\mathcal{E}^{H}}\frac{\gamma}{h}\int_{E}\overline{\kappa}\jump{u} \, \jump{v} \leq\|u\|_{\text{DG}}\,\|v\|_{\text{DG}}.
\]
For an interior coarse edge $E\in\mathcal{E}^{H}$, we let $K^{+},K^{-}\in\mathcal{T}^{H}$
be the two coarse grid blocks having the edge $E$. By the Cauchy-Schwarz
inequality, we have
\begin{equation}
\int_{E}\average{{\kappa} \nabla u\cdot{n}_{E}}\cdot\jump{v}\leq\Big(h\int_{E}\average{{\kappa}\nabla{u}\cdot{n}_{E}}^{2}(\overline{\kappa})^{-1}\Big)^{\frac{1}{2}}\Big(\frac{1}{h}\int_{E}\overline{\kappa}\jump{v}^{2}\Big)^{\frac{1}{2}}.\label{eq:cont1}
\end{equation}
Notice that
\begin{equation*}
 h\int_{E}\average{{\kappa}\nabla{u}\cdot{n}_{E}}^{2}(\overline{\kappa})^{-1}
\leq  h\Big(\int_{E}({\kappa}^{+}\nabla{u}^{+}\cdot{n}_{E})^{2}(\kappa_{K^+})^{-1}+\int_{E}({\kappa}^{-}\nabla{u}^{-}\cdot{n}_{E})^{2}(\kappa_{K^-})^{-1}\Big)
\end{equation*}
where ${u}^{\pm}={u}|_{K^{\pm}}$, $\kappa^{\pm}=\kappa|_{K^{\pm}}$.
So, summing the above over all $E$, we have
\begin{equation*}
h\sum_{E\in\mathcal{E}^{H}}\int_{E}\average{{\kappa}\nabla{u}\cdot{n}_{E}}^{2}(\overline{\kappa})^{-1}  \leq h\sum_{K\in\mathcal{T}_{H}}\int_{\partial K}(\kappa \nabla{u}\cdot{n}_{\partial K})^{2} (\kappa_{K})^{-1}
  \leq  C_{\text{inv}}^{2}a_{H}(u,u).
\end{equation*}
Thus we have
\begin{equation}
\sum_{E\in\mathcal{E}^{H}}\int_{E}\average{{\kappa}\nabla{u}\cdot{n}_{E}}\jump{v}\leq C_{\text{inv}}a_{H}(u,u)^{\frac{1}{2}}\Big(\sum_{E\in\mathcal{E}^{H}}\frac{1}{h}\int_{E}\overline{\kappa}\jump{v}^{2}\; ds\Big)^{\frac{1}{2}}.
\label{eq:cont5}
\end{equation}
Similarly, we have
\[
\sum_{E\in\mathcal{E}^{H}}\int_{E}\average{{\kappa}\nabla{v}\cdot{n}_{E}}\jump{u} \leq C_{\text{inv}}a_{H}(v,v)^{\frac{1}{2}}\Big(\sum_{E\in\mathcal{E}^{H}}\frac{1}{h}\int_{E}\overline{\kappa}\jump{u}^{2}\; ds\Big)^{\frac{1}{2}}.
\]
Summing the above two inequalities, we have
\begin{equation}
\sum_{E\in\mathcal{E}^{H}}\int_{E}\Big(\average{{\kappa}\nabla{u}\cdot{n}_{E}}\jump{v}+\average{{\kappa}\nabla{v}\cdot{n}_{E}}\jump{u}\Big)\leq C_{\text{inv}}\gamma^{-\frac{1}{2}}\|u\|_{\text{DG}}\,\|v\|_{\text{DG}}.\label{eq:cont4}
\end{equation}
This proves the continuity (\ref{lem:cont}).

For the coercivity (\ref{lem:coer}), we have
\[
a_{\text{DG}}({u},{u})=\|u\|_{\text{DG}}^{2}-\sum_{E\in\mathcal{E}^{H}}\int_{E}\Big(\average{\kappa\nabla u\cdot{n}_{E}}\cdot\jump{u}+\average{\kappa\nabla u\cdot{n}_{E}}\cdot\jump{u}\Big).
\]
By (\ref{eq:cont4}), we have
\[
a_{\text{DG}}({u},{u})\geq(1-C_{\text{inv}}\gamma^{-\frac{1}{2}})\|u\|_{\text{DG}}^{2},
\]
which gives the desired result.

\bibliographystyle{plain}
\bibliography{references}
\end{document}